\documentclass[12pt,twoside]{article}%
\usepackage{enumerate}
\usepackage{amsmath,amssymb,amsfonts}
\usepackage{amsthm}%
\usepackage{amsmath}%
\usepackage{amsfonts}%
\usepackage{amssymb}%
\usepackage{graphicx}
\usepackage{color}
\providecommand{\U}[1]{\protect\rule{.1in}{.1in}}
\newtheorem{thm}{Theorem}[section]
\newtheorem{lem}[thm]{Lemma}

\newtheorem{prop}[thm]{Proposition}
\newtheorem{corollary}[thm]{Corollary}

\newcommand{\fim}{\hfill\rule{2mm}{2mm}}
\numberwithin{equation}{section}

\newcommand{\dist}{\operatorname{dist}}
   \newcommand{\cqd}{\hfill $\blacksquare$ \vspace{0.5cm}}
   \newcommand{\rn} {\mathbb{R}^N}
   \newcommand{\dd}{\,\mathrm{d}}
   \newcommand{\cinfsc}{C_c^\infty}
   \newcommand{\grad}{\nabla}
   \newcommand{\lap}{\Delta}
   \newcommand{\ol}{{\Omega_\lambda}}
   
   \newcommand{\al}{{A_\lambda}}

   \newcommand{\imerso}{\hookrightarrow}
   \newcommand{\tofraco}{\rightharpoonup}

   \newcommand{\re}{\mathbb R}

\begin{document}
\title{The use of the Morse theory to estimate the number of nontrivial solutions of a nonlinear Schr\"odinger equation with magnetic field}
\author{Claudianor O. Alves\footnote{C.O. Alves was partially supported by CNPq/Brazil  301807/2013-2. Email: coalves@dme.ufcg.edu.br},\quad  Rodrigo C. M. Nemer \thanks{R.C.M. Nemer was supported by   FAPESP 2010/05892-9 and
CAPES/PNPD- UFCG/Matem\'atica. E-mail:rodrigocmnemer@gmail.com}\\
\noindent Unidade Acad\^emica de Matem\'atica \\
\noindent Universidade Federal de Campina Grande \\
\noindent 58429-900, Campina Grande - PB, Brazil.\\
\mbox{} and
\mbox{}\\
\noindent S\'{e}rgio H. M. Soares\footnote{Corresponding author.
 E-mail: {monari@icmc.usp.br}, Phone: +55\,16\,3373\,9660, Fax: +55\,16\,3373\,9650.}\\
\noindent Departamento de Matem\'atica \\
\noindent Instituto de Ci\^{e}ncias Matem\'{a}ticas e de
 Computa\c{c}\~{a}o\\
\noindent Universidade de S\~ao Paulo \\
\noindent 13560-970, S\~ao Carlos - SP, Brazil.
 }
\date{}
\maketitle

\begin{abstract}
Based on some ideas introduced by Benci and Cerami \cite{benci-cerami_calvar}, we obtain an abstract result that establishes a version of the Morse relations. Afterward, we use this result to prove multiplicity of solutions for a nonlinear Schr\"odinger equation with an external magnetic field.
\end{abstract}



\section{Introduction}

The relations between topological properties of the domain and the number of solutions of elliptic problems have been extensively studied by many authors. In 1991, Benci and Cerami in the pioneer paper \cite{benci-cerami_arma} studied the existence and multiplicity of solutions for the problem
\begin{equation}
\left\{\begin{array}{rcl}
            -\Delta u +  \kappa u & = & |u|^{p-2}u \quad \textrm{in } \Omega,\\
            u & = & 0 \quad \textrm{on } \partial\Omega,
       \end{array}
\right.\label{P1}
\end{equation}
where $\kappa\in \mathbb{R}^+\cup \{0\}$,  $\Omega \subset \mathbb{R}^{N}$ is a bounded smooth domain, $p \in (2,2^{*})$ and $2^*={2N}/{(N-2)}$ with $N \geq 3$. It was proved that  (\ref{P1}) has at least ${\rm{cat}}(\Omega)$ positive solutions provided  that  $\kappa$ is sufficiently large or $p$ is sufficiently close to $2^*$, where ${\rm{cat}}(\Omega)$ denotes the   Ljusternik-Schnirelman category of $\overline{\Omega}$  in itself.

Subsequently, in 1994,  Benci and Cerami in \cite{benci-cerami_calvar}   showed that the number of positive solutions for a semilinear elliptic equations like
\begin{equation}
\left\{\begin{array}{rcl}
            -\varepsilon\Delta u +  u & = & f(u) \quad \textrm{in } \Omega,\\
            u & = & 0 \quad \textrm{on } \partial\Omega,
       \end{array}
\right. \label{P2}
\end{equation}
where $\varepsilon\in  \mathbb{R}^+\setminus \{0\}$,   $\Omega \subset \mathbb{R}^{N}$ is a bounded smooth domain and $f$ is a continuous function with subcritical growth, depends on the Poincar\'{e} polynomial of the domain, that is, a lower estimate of the number of solutions can be performed entirely in terms of the Morse relations. More precisely, the authors proved   among other things that there exists $\varepsilon^* >0$ such that, for any $\varepsilon \in (0, \varepsilon^*)$ problem  (\ref{P2})  has at least $2 \mathcal P_1(\Omega) - 1$ nontrivial solutions, where $\mathcal{P}_t(\Omega)$ denotes the Poincar\'e polynomial of $\Omega$.

Multiplicity of solutions by the use of  Ljusternik-Schnirelman category or
Morse theory has been considered for different classes of problems by several authors since the works   \cite{benci-cerami_arma, benci-cerami_calvar},  see for example, Benci \cite{benci}, Benci, Bonanno and  Micheletti \cite{benci-bonanno-micheletti}, Cerami and Wei \cite{cerami-wei}, Cingolani \cite{congolani2}, Cingolani and Clapp \cite{cingolani-clapp},  Clapp \cite{clapp}, Furtado \cite{furtado}, Ghimenti and Micheletti \cite{ghimenti-micheletti}, He \cite{he}, Shang and Zhang \cite{shang-zhang} and their references.

The present paper was mainly motivated by  \cite{benci-cerami_calvar}.
By carefully examining the method used by Benci and Cerami to study some properties of the functional associated with (\ref{P2}) to apply the Morse relations,
we have observed there is an abstract result behind this method providing these relations and which can be proved by adapting the argument employed in that paper. To illustrate, we apply this result to estimate the number of nontrivial solutions for a nonlinear Schr\"odinger equations with an external magnetic field. We believe that this abstract result can be useful for finding solutions for a wide variety of elliptic problems.


In order to establish the abstract result, we need to fix some notations. Let $(E, \langle ~,~ \rangle)$ denote a real Hilbert space endowed with the induced norm $\|\cdot\|^2 = \langle ~,~ \rangle$. Let $I: E \to \re$ be a $C^2$ functional and let ${\mathcal M}$ be the Nehari manifold associated with
  $I$ given by
$$
{\mathcal M} = \{u \in E\backslash\{0\} ; I'(u)u = 0\}.
$$
Here $I$ is assumed to be bounded from below on ${\mathcal M}$ and set
\begin{equation}\label{b}
b = \inf_{\mathcal M} I.
 \end{equation}
For $a \in \mathbb{R}$, consider the sets
$$
I^{a} = \{u \in E ; I(u) \le a\} ~~~~ \mbox{and} ~~~~ {\mathcal M}^a = {\mathcal M} \cap I^a.
$$
We can now state the above-mentioned abstract result.

\begin{thm} \label{resabs}
For $b$ given by (\ref{b}), let $\delta \in (0, b)$. Suppose that
\begin{itemize}
\item[$(i)$] $I(u) = \frac{1}{2}\|u\|^2 - \Psi(u)$, where $\Psi:E \to \re$ is such that $\Psi(0) = 0$ and $t  \mapsto \Psi'(tu)u/t$ is strictly increasing in  $(0,    +\infty)$ and unbounded above, for every $u \in E\backslash \{0\}$,

\item[$(ii)$] $I$ satisfies the Palais-Smale condition and, for every $u \in E$, there exists a self-adjoint operator $L(u): E\to E$ such that $H_I(u)(v, v) = \langle L(u) v, v\rangle_E$, for every $v \in E$, where $H_I$ is the Hessian form of I at $u$,
\item[$(iii)$] The Nehari manifold ${\mathcal M}$ is homeomorphic  to the unit sphere in $E$,
\item[$(iv)$] There exist a regular value $b^* >b$ of $I$, a nonempty set $\Theta \subset \rn$ with smooth boundary and continuous applications $\Phi: \Theta^- \to {\mathcal M}^{b^*}$,  $\beta: {\mathcal M}^{b^*} \to \Theta^+$ such that $\beta \circ \Phi = Id_{\Theta^-}$, where  $$
\Theta^+ = \{ x \in \rn ; \dist (x, \Theta) \le r \}\ \mbox{ and  }\ \Theta^- = \{ x \in \Theta ; \dist (x, \partial\Theta) \ge r\},$$
for some $r > 0$ such that $\Theta^+$ and $\Theta^-$ are homotopically equivalent to $\Theta$.
\end{itemize}
Suppose also that the set $\mathcal K$ of critical points of $I$ is discrete. Then
\begin{equation}
\label{7.8}
\sum_{u \in \, \mathcal C_1} i_t (u) = t \mathcal P_t(\Theta) + t \mathcal Q(t) + (1+t)\mathcal Q_1(t)
\end{equation}
and
\begin{equation}
\label{7.9}
\sum_{u \in \, \mathcal C_2} i_t (u) = t^2 [\mathcal P_t(\Theta) + \mathcal Q(t) -1 ] + (1+t)\mathcal Q_2(t),
\end{equation}
where  $i_t (u)$ is the polynomial Morse index of $u$,
\[
\mathcal C_1 := \{ u \in \mathcal K; \delta < I(u) \le b^* \}, ~~
\mathcal C_2 := \{ u \in \mathcal K; b^* < I(u)\},
\]
$\mathcal P_t (\Theta)$ is the Poincar\'e polynomial of $\Theta$ and $\mathcal Q, \mathcal Q_1, \mathcal Q_2$ are polynomials with non-negative coefficients.
\end{thm}

As an example of the use of this result, we consider a class of nonlinear Schr\"odinger equations with an external magnetic field, namely
\begin{equation}
\left\{\begin{array}{rcl}
            \left( \frac{1}{i}\nabla - A \right)^2 u + \kappa u & = & |u|^{p-2} u \quad \textrm{in } \Omega,\\
            u & = & 0 \quad \textrm{on } \partial\Omega,
       \end{array}
\right. \label{akp}
\end{equation}
where $\kappa$ is a positive parameter, $\Omega \subset \mathbb R ^N$ is a smooth bounded domain, $N \geq 3$, $i$ is the imaginary unit and $p \in (2, 2^*)$, $2^* = 2N/(N-2)$. The function $A: \Omega \to \rn$ is the magnetic potential and the Schr\"odinger operator is defined by
\[
 \left( \frac{1}{i}\nabla - A \right)^2 u = -\Delta u - \frac{2}{i}A\cdot \nabla u + |A|^2u - \frac{1}{i}u{\rm{div}}\, A.
\]
We assume that $A \in L^\infty(\Omega, \rn)$.

Existence results for the magnetic case, that is $A \not=0$, has also received a special attention in the last year. Associated with this subject, the reader can find interesting results in the papers \cite{abatangelo-terracini},
\cite{alves-nemer-soares},
\cite{alves-figueiredo-furtado},
\cite{barile},
\cite{cao-tang},
\cite{chabrowski-szulkin},
\cite{cingolani},
\cite{cingolani-clapp}
\cite{cingolani-jeanjean-secchi1},
\cite{cingolani-jeanjean-secchi2},
\cite{cingolani-secchi}
\cite{esteban-lions},
 \cite{kurata},
\cite{li-peng-wang},
\cite{liang-zhang},
\cite{squassina},
\cite{tang1},
\cite{tang2},
\cite{tang3}.

Motivated by \cite{benci-cerami_arma, benci-cerami_calvar},   we obtain the following result.
\begin{thm} \label{t1} Suppose that the set $\mathcal K$ of solutions of the problem \eqref{akp} is discrete. Then there is a function $\overline p : [ \left.0, +\infty )\right. \to (2, 2^*)$ such that for every $p \in [\left. \overline p(\kappa), 2^* )\right.$,
\[
\sum\limits_{u \in \mathcal K} i_t(u)= t \mathcal P_t (\Omega) + t^2 [\mathcal P_t (\Omega) - 1] + 
\mathcal Q (t),
\]
where $\mathcal Q$ is a polynomial with non-negative integer coefficients, $\mathcal P_t (\Omega)$ is the Poincar\'{e} polynomial of $\Omega$ and $i_t(u)$ is the Morse index of $u$.
\end{thm}

In the non-degenerate case, we have:

\begin{corollary} \label{c1}
Suppose that the solutions of problem \eqref{akp} are non-degenerate. Then there is a function $\overline p : [ \left.0, +\infty )\right. \to (2, 2^*)$ such that for every $p \in [\left. \overline p(\kappa), 2^* )\right.$, problem \eqref{akp}  has at least $2\mathcal P_1 (\Omega) -1$ nontrivial solutions.
\end{corollary}

Another application of the abstract result can be given by the following problem
\begin{equation}
\left\{\begin{array}{rcl}
            \left( -i\nabla - A_\lambda \right)^2 u + u & = & |u|^{p-2}u, \quad \textrm{in } \Omega_\lambda, \\
            u & = & 0, \quad \textrm{on } \partial\Omega_\lambda,
       \end{array}
\right.\label{pal0}
\end{equation}
where $\lambda > 0$ is a positive parameter,   $A_\lambda := A\left({x}/{\lambda}\right)$, $\Omega_\lambda := \lambda \Omega$, $\Omega \subset \mathbb R ^N$ $(N \ge 3)$  is a bounded smooth domain and $p \in (2, 2^*)$ .  We observe that, unlike the case with no magnetic vector field $A$, problem (\ref{pal0}) cannot be written in the form (\ref{akp}), and hence these problems are different.  In \cite{alves-figueiredo-furtado}, Alves et al have proved that for large values of $\lambda>0$, problem (\ref{pal0}) has at least ${\rm{cat}}(\Omega_\lambda)$ nontrivial weak solutions. Combing the abstract result with arguments present in \cite{alves-figueiredo-furtado}, we are able to estimate the number of nontrivial solution in terms of the $\mathcal P_t(\Omega_{\lambda})$. More precisely, we can prove that  (\ref{pal0})   has at least $2 \mathcal P_1(\Omega_{\lambda}) - 1$ nontrivial solutions provided that  $\lambda$ is sufficiently large.

\section{The abstract theorem }

In this section we give the proof of Theorem \ref{resabs}. We begin by showing how the set $\Theta$ relates to the set ${\mathcal M}^{b^*}$.

\begin{lem} \label{al1} Under the assumptions of Theorem \ref{resabs}, we have
\[
\mathcal P _t ({\mathcal M}^{b^*}) = \mathcal P_t (\Theta) + \mathcal Q(t),
\]
where $\mathcal Q$ is a polynomial with non-negative coefficients.
\end{lem}

\noindent {\bf Proof.}~ We observe that $\Phi$ induces a homomorphism $(\Phi)_k: H_k(\Theta^-) \to H_k({\mathcal M}^{b^*})$ between the $k$-th homology groups. Since $\Phi$  is a {{injective}} function, so also is  $(\Phi)_k$. Hence, $\operatorname{dim}H_k(\Theta^-) \geq \operatorname{dim}H_k({\mathcal M}^{b^*})$, and the result follows from the definition of the Poincar\'{e} polynomials and the fact that $\Theta^-$ and $\Theta$ are homotopically equivalent.  \fim

\begin{lem} \label{al2} Let $\delta \in (0, b)$ and let  $a \in \left.(\delta, \infty\right.]$ be a noncritical level of $I$. Then
\[
\mathcal P_t (I^a, I^\delta) = t \mathcal P_t({\mathcal M}^a).
\]
\end{lem}

\noindent {\bf Proof.} The proof proceeds along the same lines as the proof of \cite[Lemma 5.2]{benci-cerami_calvar}.\fim

\begin{lem} \label{al8} Let $\delta$ be as in Lemma \ref{al2}. Then
\begin{equation}
\mathcal P_t (I^{b^*}, I^\delta) = t \mathcal P_t (\Theta) + t \mathcal Q (t)
\label{al81}
\end{equation}
and
\begin{equation}
\mathcal P_t (E, I^\delta) = t \mathcal P_t ({\mathcal M})=t,
\label{al82}
\end{equation}
where $\mathcal Q$ is a polynomial with non-negative coefficients.
\end{lem}

\noindent {\bf Proof.}  By assumption, $b^*$ is a regular value. Applying Lemma \ref{al2}, for  $a = b^*$, and Lemma \ref{al1}, we get \eqref{al81}. Using the fact that ${\mathcal M}$ is homeomorphic to the unit sphere in $E$, which we know to be contractible, we have that ${\mathcal M}$ is contractible. Hence, $\operatorname{dim}H^k({\mathcal M}) = 1$ if $k = 0$ and $\operatorname{dim}H^k({\mathcal M}) = 0$ if $k \neq  0$. The identity
 \eqref{al82} follows from Lemma \ref{al2} with  $a = +\infty$ and the fact that ${\mathcal M}$ is contractible. \fim

\begin{lem} \label{al9}
Let $\delta$ be as in Lemma \ref{al2}. Then
\begin{equation}
\mathcal P_t (E, I^{b^*}) = t^2 [ \mathcal P_t (\Theta) + \mathcal Q(t) - 1 ], \label{al91}
\end{equation}
where $\mathcal Q$ is a polynomial with non-negative coefficients.
\end{lem}

\noindent {\bf Proof.} We follow Benci and Cerami \cite{benci-cerami_calvar} in considering the exact sequence:
\[
\dots \longrightarrow H_k (E, I^\delta) \stackrel{j_k}\to  H_k (E, I^{b^*})
\stackrel{\partial_k}\to H_{k-1} (I^{b^*}, I^\delta) \stackrel{i_{k-1}}\to H_{k-1} (E , I^\delta) \longrightarrow \dots
\]
From \eqref{al82}, we obtain $\operatorname{dim} H_k(E, I^\delta) = 0$, for every $k \neq 1$. Combining this with the fact that the sequence is exact, we obtain that  $\partial_k$ is an isomorphism for every $k \ge 3$. Hence,
\begin{equation}
\operatorname{dim} H_k(E, I^{b^*}) = \operatorname{dim} H_{k-1} (I^{b^*}, I^\delta), \, \forall k \ge 3.
\label{al92}
\end{equation}
For $k = 2$, we have
\[
\dots \longrightarrow H_2 (E, I^\delta) \stackrel{j_2}\to  H_2 (E, I^{b^*}) \stackrel{ \partial_2}\to H_{1} (I^{b^*}, I^\delta) \stackrel{i_1}\to H_{1} (E, I^\delta) \longrightarrow \dots
\]
Since the homomorphism induced by the canonic projection $j_2$  is surjective and $\operatorname{dim} H_2 (E, I^\delta) = 0$, by \eqref{al82}, we have
\begin{equation}
H_2 (E, I^{b^*}) = j_2 (H_2(E, I^\delta)) = \{0\}.
\label{al93}
\end{equation}
For $k = 1$,
\[
\dots \longrightarrow H_1 (I^{b^*}, I^\delta) \stackrel{i_1}\to H_1 (E, I^\delta) \stackrel{j_1}\to H_1 (E, I^{b^*}) \stackrel{\partial_1}\to H_0(I^{b^*}, I^\delta) \longrightarrow \dots
\]
Using that $E$ is a connected set, we have
\begin{equation}
H_0 (E, I^{b^*}) = 0.
\label{al94}
\end{equation}
We now claim that $i_1$ is an isomorphism. Indeed, as $\Theta \neq \emptyset$ and  $\operatorname{dim} H_0 (\Theta)$ is the number of connected components of the set  $\Theta$, we have $H_0 (\Theta) \neq \{0\}$. By \eqref{al81}, $H_1 (I^{b^*}, I^\delta) \neq \{0\}$. From \eqref{al82}, we obtain $\operatorname{dim} H_1 (E, I^\delta) = 1$. Since $i_1$ is injective, it follows that $\operatorname{dim} H_1 (I^{b^*}, I^\delta) = 1$, and so $i_1$ is an isomorphism. Hence, as $j_1$ is surjective, we get
\begin{equation}
\operatorname{dim} H_1 (E, I^{b^*}) = 0.
\label{al95}
\end{equation}
Combining Lemma \ref{al8} with \eqref{al92}-\eqref{al95}, we have
\begin{eqnarray*}
\mathcal P_t (E, I^{b^*}) & = & \sum_{k \ge 3} t^k \operatorname{dim} H_k (E, I^{b^*})\\
                          & = & \sum_{k \ge 3} t^k \operatorname{dim} H_{k-1} (I^{b^*}, I^\delta) = t\sum_{k \ge 3} t^{k-1} \operatorname{dim} H_{k-1} (I^{b^*}, I^\delta) \\
                          & = & t \left[ \mathcal P_t (I^{b^*}, I^\delta) - t \, \operatorname{dim}H_1 (I^{b^*}, I^\delta) - \operatorname{dim}H_0 (I^{b^*}, I^\delta) \right]\\
                          & = & \, t^2 \left[ \mathcal P_t (\Theta) + \mathcal Q(t) - 1 \right],
\end{eqnarray*}
which completes the proof of Lemma \ref{al9}. \fim

\vspace{0.5 cm}

\subsection{Proof of Theorem \ref{resabs}}

Now, we are able to conclude proof of Theorem \ref{resabs}.  By $(ii)$, $I$ satisfies the Palais-Smale condition and, for a nondegenerate critical point $u$ of $I$, the linear operator $L(u)$ associated to $H_I(u)$ is a Fredholm operator with index 0.  By \cite[Example I.5.1]{benci}, we can use \cite[Theorem I.5.9]{benci} and Lemmas \ref{al8} and \ref{al9} to get
\begin{eqnarray*}
\sum_{u \in \mathcal C _1} i_t (u) & = & \mathcal P _t (I^{b^*}, I^\delta) + (1+t)\mathcal Q_1(t)\\
                                   & = & t \big[\mathcal P_t (\Theta) + \mathcal Q (t)\big] + (1+t)\mathcal Q_1(t)
\end{eqnarray*}
and
\begin{eqnarray*}
\sum_{u \in \mathcal C _2} i_t (u) & = & \mathcal P _t (E, I^{b^*}) + (1+t)\mathcal Q_2(t)\\
                                   & = & t^2 \big[\mathcal P_t (\Theta) + \mathcal Q (t) - 1 \big] + (1+t)\mathcal Q_2(t).
\end{eqnarray*} \fim


\section{Application of the abstract theorem}

This section is devoted to prove Theorem \ref{t1}.  Let $E$ be a real Hilbert space defined as the closure of $\cinfsc(\Omega, \mathbb C)$ with respect to the norm induced by the inner product
\[
\langle u, v\rangle_\kappa := \textrm{Re} \left\{ \int_{\Omega} \big[ \grad_A u \overline{\grad_A v} + \kappa u \overline v \big] d x \right\},
\]
where, for $a, b \in \mathbb C^M$, $M \in \mathbb N$, $ab = \sum_{j=1}^M a^j.b^j$, where ``\,.\,'' is the usual complex multiplication, $\textrm{Re} (a)$ is the real part of $a \in \mathbb C^M$ and $\overline a$ the complex conjugate of $a$. Moreover,
\[
\grad_A u := (D^j_{A} u)_{j = 1}^{N}, \quad D^j_{A} u = -i\partial_j u - A^j u, \quad j \in {1, ..., N}.
\]
The norm induced by this inner product is
\[
\|u\|_\kappa^2 := \int_\Omega \big[ |\grad_A u|^2 + \kappa |u|^2 \big] d x.
\]
As proved in Esteban and Lions \cite{esteban-lions}, for every $u \in E$ there holds
\[
|\grad_A u| \ge |\grad |u||.
\]
The above expression is the so called diamagnetic inequality. The functional associated with \eqref{akp}, $I_{\kappa, p, \Omega}: E \to \mathbb R$, is given by
\[
I_{\kappa, p, \Omega}(u) = \frac{1}{2} \int_{\Omega} (|\grad_A u|^2 + \kappa |u|^2) d x- \frac{1}{p} \int_{\Omega} |u|^p d x, \forall u \in E.
\]
By Sobolev embeddings and diamagnetic inequality, $I_{\kappa, p, \Omega}$ is well defined. Furthermore, $I_{\kappa, p, \Omega} \in C^2 (E, \mathbb R)$ with
\[
I_{\kappa, p, \Omega}'(u)v = \mathrm{Re}\left( \int_{\Omega} (\grad_A u \overline{\grad_A v} + \kappa u\overline v - |u|^{p-2}u\overline v) d x \right), \forall u, v \in E.
\]
Thus, every critical point of $I_{\kappa, p, \Omega}$ is a weak solution of \eqref{akp}.

\vspace{0.5 cm}

A standard verification shows that:

{\prop \label{ps} The functional $I_{\kappa, p, \Omega}$ satisfies Palais-Smale condition, that is, every sequence $(u_n) \subset E$ for which $\sup_n |I_{\kappa, p, \Omega}(u_n)| < +\infty$ and $I_{\kappa, p, \Omega}'(u_n) \to p$, as $n \to \infty$, has a convergent subsequence.}

\medskip
It is straightforward to show that $I_{\kappa, p, \Omega}$ satisfies the geometric hypotheses of the mountain pass theorem. From this and Proposition \ref{ps}, for all $p \in (2, 2^*)$ and $\kappa > 0$, problem \eqref{akp} has a nontrivial solution $u \in E$ such that $I_{\kappa, p, \Omega}(u) = b_{\kappa, p, \Omega}$ and $I_{\kappa, p, \Omega}' (u) = 0$, where $b_{\kappa, p, \Omega}$ denotes the mountain pass level $I_{\kappa, p, \Omega}$. Moreover, as in \cite[Theorem 4.2]{Willem},
\[
b_{\kappa, p, \Omega} := \inf_{u \in {\mathcal M}_{\kappa, p, \Omega}} I_{\kappa, p, \Omega} (u),
\]
where
$
{\mathcal M}_{\kappa, p, \Omega} = \left\{ u \in E \backslash \{0\}; I_{\kappa, p, \Omega}' (u)u = 0\right\}
$
is the Nehari manifold associated to $I_{\kappa, p, \Omega}$.

{\prop \label{nehari} The Nehari manifold ${\mathcal M}_{\kappa, p, \Omega}$ is diffeomorphic  to the unit sphere of $E$. Moreover, there is $\delta=\delta(p)> 0$, independent of $\kappa > 0$, such that for every $u \in {\mathcal M}_{\kappa, p, \Omega}$,
\begin{equation*}
\int_{\Omega} (|\grad_A u|^2 + \kappa |u|^2) d x \ge \delta \quad \textrm{and} \quad I_{\kappa, p, \Omega}(u) \ge \delta.
\end{equation*}}

\noindent {\bf Proof.} For any $u \in {\mathcal M}_{\kappa, p, \Omega}$, the diamagnetic inequality combined with Sobolev imbedding imply
\begin{align*}
\|u\|^2 & = \int_{\Omega} |u|^p d x \le C_p \int_{\Omega} |\grad |u||^2 \le C_p\|u\|^p,
\end{align*}
where $C_p$ is the constant of the embedding $H^1_0(\Omega) \imerso L^p(\Omega, \mathbb R)$. Thus
\[
\|u\| \ge {C_p}^{\frac{1}{2-p}} =: \delta_1,
\]
from where it follows
\[
I_{\kappa, p, \Omega} (u) = \left(\frac{1}{2} - \frac{1}{p}\right)\|u\|^2
                 \ge \left(\frac{1}{2} - \frac{1}{\theta}\right) \delta_1^2 =: \delta.
\]
To conclude the proof, let $S$ be the unit sphere in $E$. For every $u \in S$,  let
$\xi(u) > 0$  be the unique positive number such that
\[
\frac{d }{d t} I_{\kappa, p, \Omega}(tu)\big|_{t=\xi(u)} = 0.
\]
This define a $C^1$ function  $\xi : S \to (0, +\infty)$  by the Implicit Function Theorem.  Thus, $D: S \to {\mathcal M}_{\kappa, p, \Omega}$ given by
\[
D(u) = \xi(u)u \in {\mathcal M}_{\kappa, p, \Omega}
\]
is a $C^1$ diffeomorphism. \fim

{\prop \label{psr} $I_{\kappa, p, \Omega} \big|_{{\mathcal M}_{\kappa, p, \Omega}}$ satisfies the Palais-Smale condition}.

\noindent {\bf Proof.} Let $(u_n) \subset {\mathcal M}_{\kappa, p, \Omega}$ be a sequence satisfying
$$
\sup_{n \in \mathbb N} |I_{\kappa, p, \Omega}(u_n)| < \infty ~~ \mbox{and} ~~  \left(I_{\kappa, p, \Omega} \big|_{{\mathcal M}_{\kappa, p, \Omega}}\right)'(u_n) \to 0.
$$
Taking a subsequence if necessary, we can assume that $I_{\kappa, p, \Omega}(u_n) \to d$ as $n \to \infty$. A standard verification shows that  $(u_n)\subset E$ is bounded. Thus there is $u \in E$ such that $u_n \tofraco u$ in $E$. Consequently, $u_n \to u$ in $L^p(\Omega, \mathbb C)$. By \cite[Proposition 5.12]{Willem}, for each $n \in \mathbb N$, there is $\mu_n \in \mathbb R$ such that
\begin{equation}
\label{psr1}
I_{\kappa, p, \Omega}'(u_n) - \mu_n G'(u_n) = \left(I_{\kappa, p, \Omega} \big|_{{\mathcal M}_{\kappa, p, \Omega}}\right)'(u_n) = o_n(1),
\end{equation}
where $G_{\kappa, p, \Omega}(u) = I_{\kappa, p, \Omega}'(u)u$. Since $u_n \in {\mathcal M}_{\kappa, p, \Omega}$, by Proposition \ref{nehari}, we obtain
\begin{eqnarray*}
\lim_n G_{\kappa, p, \Omega}'(u_n)u_n = \lim_n (2 - p)\int_{\Omega}|u_n|^p dx  \le (2-p)\delta<0.
\end{eqnarray*}
This and \eqref{psr1} imply that $\mu_n \to 0$ as $n \to 0$. The result follows from Proposition \ref{ps}. \cqd

\begin{corollary} \label{pcr} If $u$ is a critical point of $I_{\kappa, p, \Omega}$ constrained to ${\mathcal M}_{\kappa, p, \Omega}$, then $u$ is a critical point of $I_{\kappa, p, \Omega}$.
\end{corollary}

\noindent {\bf Proof.} The proof proceeds along the same lines as the proof of Proposition \ref{psr}. \cqd



\subsection{Behaviour of the minimax levels }
For any $p\in (2, 2^*)$ and $\kappa >0$ we denote
\[
m_A(\kappa, p, \Omega) := \inf\limits_{u \in E \backslash \{0\}} \frac{\displaystyle\int_{\Omega} (|\grad_A u|^2 + \kappa |u|^2) \dd x }{\left(\displaystyle\int_{\Omega}|u|^p \dd x \right)^{\frac{2}{p}}},
\]
\[
S_{A, \kappa} := m_A(\kappa, 2^*, \Omega), \quad \mbox{and}\quad S_A := m_A(0, 2^*, \Omega).
\]
Employing the same arguments in \cite{Willem}, we can prove the following result:

\begin{lem} \label{bma} Let $b_{\kappa, p, \Omega}$ be the mountains pass level of $I_{\kappa, p, \Omega}$. Then
\[
b_{\kappa, p, \Omega} = \left(\frac{1}{2} - \frac{1}{p}\right) m_A(\kappa, p, \Omega)^{\frac{p}{p-2}}.
\]
Hence,
\[
b_{\kappa, 2^*, \Omega} = \frac{1}{N} {S_{A, \kappa}}^{\frac{N}{2}}.
\]
\end{lem}
From now on, we also consider
\[
m(\kappa, p, \Omega) := \inf\limits_{u \in H^1_0(\Omega) \backslash \{0\}} \frac{\displaystyle\int_{\Omega}(|\grad u|^2 + \kappa u^2) \dd x }{\left(\displaystyle\int_{\Omega}|u|^p \dd x \right)^{\frac{2}{p}}}
\]
and
\[
S_\kappa := m(\kappa, 2^*, \Omega).
\]
Then, $S := m(0, 2^*, \Omega)$, where $S$ is the best constant of the imbedding $H_0^{1}(\Omega, \re) \hookrightarrow L^{2^*}(\Omega, \re)$, which is independent of $\Omega$. Moreover, from \cite[Theorem 1.1]{arioli-szulkin}, we have

\begin{lem} \label{ssal} For every $\kappa \ge 0$, we have $S_{A, \kappa} = S_\kappa = S$.
\end{lem}

The following lemma is the key to establish a  relation between $b_{\kappa, p, \mathcal D} $ and $b_{2^*}$.

\begin{lem}  \label{bkpd} For any given $\kappa \ge 0$ and for any  bounded domain $\mathcal D \subset \rn$, the following limit holds:
\[
\lim\limits_{p \to 2^*} b_{\kappa, p, \mathcal D} = b_{2^*},
\]
where $b_{2^*}$ denotes the mountain pass level associated with the functional $J_{\infty}:H_0^{1}(\Omega) \to \mathbb{R}$ given by
\[
J_\infty(u) = \frac{1}{2} \int_{\Omega} |\grad u|^2 d x- \frac{1}{2^*} \int_{\Omega} |u|^{2^*} d x.
\]

\end{lem}

\noindent {\bf Proof.} ~ Fix $\kappa \ge 0$ and $\mathcal D \subset \rn$  a bounded domain. Now let  $2 \le p < q \le 2^*$ and $u \in E(\mathcal D)$, where the Hilbert space $E(\mathcal D)$ is defined of the same way of $E$ taking $\mathcal D$ instead of $\Omega$. 
Notice that $|u|_{p, \Omega} \le |\mathcal D|^{\frac{q-p}{qp}} |u|_{q, \mathcal D}$, so
\begin{equation} \label{bkpd3}
\frac{\displaystyle \int_{\mathcal D} (|\grad_A u|^2 + \kappa |u|^2) dx}{|u|^2_{p, \mathcal D}} \ge |\mathcal D|^{\frac{-2(q-p)}{qp}} \frac{\displaystyle \int_{\mathcal D} (|\grad_A u|^2 + \kappa |u|^2) d x}{|u|^2_{q, \mathcal D}}.
\end{equation}
Taking $q = 2^*$ and the infimum over all $u \in E(\mathcal D)\backslash\{0\}$, we find
\begin{equation}
\label{bkpd1}
m_A(\kappa, p, \mathcal D) \ge |\mathcal D|^{\frac{-2(2^*-p)}{2^*p}} S_A.
\end{equation}
On the other hand, taking $p = 2$, $q = p$ and using similar arguments, we obtain
\begin{equation}
\label{bkpd2}
m_A(\kappa, p, \mathcal D) \le |\mathcal D|^{\frac{p-2}{p}} m_A(\kappa, 2, \mathcal D).
\end{equation}
Then, by \eqref{bkpd1} and \eqref{bkpd2}, $\left(m_A(\kappa, p, \mathcal D)\right)_p$ is a bounded sequence, therefore there exist
$$
M:= \limsup\limits_{p \to 2^*} m_A(\kappa, p, \mathcal D) ~~ \mbox{and} ~~ m:= \liminf\limits_{p \to 2^*}m_A(\kappa, p, \mathcal D).
$$
We claim that $M = S_A = m$. Indeed, by (\ref{bkpd1}),
$$
m \ge \liminf\limits_{p \to 2^*} |\mathcal D|^{\frac{-2(2^* - p)}{2^* p}} S_A = S_A.
$$
Suppose by contradiction that $m > S_A$.  Let $\epsilon \in (0, m - S_A)$. By the definition of $S_A$, there is $\overline u \in E(\mathcal D)$ such that
\[
\frac{\displaystyle\int_{\mathcal D} (|\grad_A \overline u|^2 + \kappa |\overline u|^2) d x }{\left(\displaystyle\int_{\mathcal D} |\overline u|^{2^*} d x \right)^{\frac{2}{2^*}}} < S_A + \frac{\epsilon}{2}.
\]
On the other hand, as the  function $p \mapsto |\overline u|_{p, \mathcal D}$ is continuous, there exists $\overline p \in (2, 2^*)$ such that for every $p \in [\left. \overline p, 2^* )\right.$, we have
$$
\left| \frac{\displaystyle\int_{\mathcal D} (|\grad_A \overline u|^2 + \kappa |\overline u|^2)d x }{\left( \displaystyle\int_{\mathcal D} |\overline u|^p d x \right)^{\frac{2}{p}}} - \frac{\displaystyle\int_{\mathcal D} (|\grad_A \overline u|^2 + \kappa |\overline u|^2) d x }{\left( \displaystyle\int_{\mathcal D} |\overline u|^{2^*} d x\right)^{\frac{2}{2^*}}} \right| < \frac{\epsilon}{2}.
$$
Thus, for every $p \in [ \overline p, 2^* ]$,
\begin{eqnarray*}
m_A(\kappa, p, \mathcal D) & \le & \frac{\displaystyle\int_{\mathcal D} (|\grad_A \overline u|^2 + \kappa |\overline u|^2 ) d x }{\left( \displaystyle\int_{\mathcal D} |\overline u|^p d x \right)^{\frac{2}{p}}} < \frac{\displaystyle\int_{\mathcal D} (|\grad_A \overline u|^2 + \kappa |\overline u|^2) d x }{\left( \displaystyle\int_{\mathcal D} |\overline u|^{2^*} d x \right)^{\frac{2}{2^*}}} + \frac{\epsilon}{2} \\
                          &  <  & S_A + \epsilon < m,
\end{eqnarray*}
that is, $m = \liminf\limits_{p \to 2^*} m_A(\kappa, p, \mathcal D) < S_A + \epsilon < m$, which is a contradiction. Hence $S_A = m$. Similar arguments show that $S_A = M$.
\cqd


In the following, for all $\kappa \ge 0$ and $p \in (2, 2^*)$, we consider the functional
\[
J_{\kappa, p, \Omega}(u) = \frac{1}{2}\int_{\Omega} (|\grad u|^2 + \kappa u^2) d x - \frac{1}{p} \int_{\Omega} |u|^p d x, \quad \forall u \in H^1_0(\Omega),
\]
and the corresponding Nehari manifold
\[
{\mathcal N}_{\kappa, p, \Omega} := \{ u \in H^1_0(\Omega)\backslash \{0\} ; J_{\kappa, p, \Omega}'(u)u = 0 \}.
\]
Define
\[
c_{\kappa, p, \Omega} = \inf\limits_{u \in {\mathcal N}_{\kappa, p, \Omega}} J_{\kappa, p, \Omega}(u).
\]
As in the proof of Lemma \ref{bma} and Proposition \ref{bkpd}, for all $\kappa \ge 0$ and $\mathcal D \subset \rn$ a bounded domain, we have
\[
c_{\kappa, p, \mathcal D} = \left(\frac{1}{2} - \frac{1}{p}\right) m(\kappa, p, \mathcal D)^{\frac{p}{p-2}} \quad \mbox{ and } \quad \lim\limits_{p \to 2^*} c_{\kappa, p, \mathcal D} = c_{\kappa, 2^*, \mathcal D}.
\]
In particular,
\[
c_{\kappa, 2^*, \mathcal D} = \frac{1}{N} {S}^{\frac{N}{2}} =: b_{2^*}.
\]
Thus, by Lemma \ref{ssal},
\begin{equation}
\label{cb}
\lim\limits_{p \to 2^*} c_{\kappa, p, \mathcal D} = \lim\limits_{p \to 2^*} b_{\kappa, p, \mathcal D} = b_{2^*}.
\end{equation}

Without loss of generality we can assume that $0 \in \Omega$. Let $r > 0$ be such that $B_r(0) \subset \Omega$ and the sets
\[
\Omega^+ := \{ x \in \mathbb R ^N; \textrm{dist}(x, \Omega) \le r \} \quad \textrm{and} \quad \Omega^- := \{ x \in \Omega; \textrm{dist}(x, \partial\Omega) \ge r \}
\]
are homotopically equivalent to $\Omega$.

Define $(I_{\kappa, p, r}; {\mathcal M}_{\kappa, p, r}; b_{\kappa, p, r})$ and $(J_{\kappa, p, r}; {\mathcal
N}_{\kappa, p, r}; c_{\kappa, p, r})$  in an exactly similar way to those of
$(I_{\kappa, p, \Omega}; {\mathcal M}_{\kappa, p, \Omega}; b_{\kappa, p, \Omega})$ and $(J_{\kappa, p, \Omega}; {\mathcal N}_{\kappa, p, \Omega}; c_{\kappa, p, \Omega})$, by taking $B_r(0) \subset \Omega$ instead of $\Omega$.

Using that $J_{\kappa, p, r}|_{{\mathcal N}_{\kappa, p, r}}$ satisfies the Palais-Smale condition,  there exists a positive function $u_{\kappa, p, r} \in {\mathcal N}_{\kappa, p, r}$ such that $J_{\kappa, p, r}(u_{\kappa, p, r})= c_{\kappa, p, r}$ and $J'_{\kappa, p, r}(u_{\kappa, p, r})=0$.
By Schwarz simmetrization we can assume that $u_{\kappa, p, r}$ is
radially symmetric.  Let  $t_{\kappa, p, y} > 0$ be the unique positive number such that   $t_{\kappa, p, y} e^{i\tau_y}u_{\kappa, p, r} (|.-y|) \in {\mathcal M}_{\kappa, p, \Omega}$.  Define the function  $\Phi_{\kappa, p}: \Omega^-_r \to {\mathcal M}_{\kappa, p, \Omega}$ as
\[
[\Phi_{\kappa, p} (y)](x) = \begin{cases}
                               t_{\kappa, p, y} e^{i\tau_y(x)}u_{\kappa, p, r} (|x-y|), & x \in B_r(y),\\
                               0, & x \in \Omega \backslash B_r(y),
                             \end{cases}
\]
where $\tau_y(x) :=  \sum^N_{j = 1} A^j(y) x^j, x = (x_1, \dots, x_N) \in \Omega$.

\begin{lem} \label{phi} For a fixed $\kappa \ge 0$,
\[
\lim\limits_{p \to 2^*} \max\limits_{y \in \Omega^-_r} |\Phi_{\kappa, p}(y) - b_{2^*}| = 0.
\]
\end{lem}
\noindent {\bf Proof.} ~Let $(p_n) \subset [2,2^*) $ and $(y_n) \subset \Omega_r^-$ be sequences such that $p_n \to 2^*$ and
\[
I_{\kappa, p_n, \Omega} (\Phi_{\kappa, p_n}(y_n)) \to b_{2^*}, \textrm{ as } n \to \infty.
\]
For simplicity, we will write
$$
t_{\kappa, p_n, y_n} =: t_n,~ I_{\kappa, p_n, \Omega} =: I_n, ~\Phi_{\kappa, p_n}(y_n) =: \Phi_n(y_n) ~~ \mbox{and} ~~ u_{\kappa, p_n, r} =: u_n.
$$
Observe that
\begin{eqnarray*}
I_n (\Phi_n(y_n)) & = & \frac{1}{2} \int_{\Omega} (|\grad_A \Phi_n(y_n)|^2 + \kappa |\Phi_n(y_n)|^2) \dd x - \frac{1}{p_n} \int_{\Omega} |\Phi_n(y_n)|^{p_n} d x\\
                                 & = & \frac{{t_n}^2}{2} \int_{B_r(0)} |A(z + y_n) - A(y_n)|^2 |u_n|^2 d x + \\
                                 & & + \frac{{t_n}^2}{2} \int_{B_r(0)} (|\grad u_n|^2 + \kappa |u_n|^2) d x - \frac{{t_n}^{p_n}}{p_n} \int_{B_r(0)} {u_n}^{p_n} d x\\
                                 & \le & \frac{{t_n}^2}{2} \int_{B_r(0)} |A(z + y_n) - A(y_n)|^2 |u_n|^2 d x+ J_{\kappa, p_n, r} (u_n)\\
                                 & = & \frac{{t_n}^2}{2} \int_{B_r(0)} |A(z + y_n) - A(y_n)|^2 |u_n|^2 d x + c_{\kappa, p_n, r}.
\end{eqnarray*}
On the other hand, by diamagnetic inequality,
\begin{eqnarray*}
I_n (\Phi_n(y_n)) & \ge & I_n( e^{i\tau_{y_n}}u_n(. - y_n) ) \\
                  & \ge & \frac{1}{2}\int_{\Omega} \big(|\grad |e^{i\tau_{y_n}}u_n(. - y_n)| |^2 + |e^{i\tau_{y_n}}u_n(. - y_n)|^2 \big) d x - \\
                  & & - \frac{1}{p} \int_{\Omega} |e^{i\tau_{y_n}}u_n(. - y_n)|^p d x \\
                  & = & J_{\kappa, p_n, r} (u_n) = c_{\kappa, p_n, r}.
\end{eqnarray*}
Thus, by \eqref{cb}, it is sufficient  to show that
\begin{equation}
\label{phi1}
\frac{{t_n}^2}{2} \int_{B_r(0)} |A(z + y_n) - A(y_n)|^2 |u_n|^2 d x = o_n(1).
\end{equation}
We begin showing that $u_n \tofraco 0$ in $H^1_0(B_r(0), \re)$ and $(t_n)_n$ is a bounded sequence. In fact, since $u_n \in {\mathcal N}_{\kappa, p_n, r}$ achieves $c_{\kappa, p_n, r}$,
\begin{equation}
\label{phi2}
\int_{B_r(0)} (|\grad u_n|^2 + \kappa |u_n|^2) d x = \left(\frac{1}{2} - \frac{1}{p_n}\right)^{-1} c_{\kappa, p_n, r}.
\end{equation}
From \eqref{cb}-\eqref{phi2}, the sequence $(u_n) \subset H^1_0(B_r(0), \re)$ is  bounded. Thus, there exists $v \in H^1_0(B_r(0))$ such that
\begin{equation}
\label{phi22}
\begin{cases}
u_n \tofraco v \mbox{ in } H^1_0(B_r(0), \re), \mbox{ as } n \to \infty\\
u_n \to v \mbox{ in } L^s(B_r(0), \re), \mbox{ for each } s \in \left.[1, 2^*\right.), \mbox{ as } n \to \infty\\
u_n(x) \to v(x) \mbox{ almost everywhere } B_r(0), \mbox{ as } n \to \infty.
\end{cases}
\end{equation}
By the fact that $u_n \in {\mathcal N}_{\kappa, p_n, r}$ achieves $c_{\kappa, p_n, r}$, $u_n$ is a solution of
\[
\begin{cases}
-\lap u + \kappa u = u^{p_n-1} \mbox{ in } B_r(0),\\
u = 0 \mbox{ on } \partial B_r(0).
\end{cases}
\]
Consequently, for any $\psi \in \cinfsc(B_r(0))$,
\[
\int_{B_r(0)} (\grad u_n \grad \psi + \kappa u_n \psi) d x = \int_{B_r(0)} {u_n}^{p-1} \psi d x.
\]
By \eqref{phi22}, as $n \to \infty$,
\begin{equation}
\label{phi3}
\int_{B_r(0)} (\grad u_n \grad \psi + \kappa u_n \psi) d x \to \int_{B_r(0)} (\grad v \grad \psi + \kappa v \psi) d x.
\end{equation}
Since $(u_n^{p_n-1})$ is a bounded sequence in $L^{\frac{2^*}{2^*-1}}(\Omega)$ and $u_n^{p_n-1}(x) \to v^{2^*-1}(x)$ almost everywhere in $\Omega$, it follows that
$$
u_n^{p_n-1} \rightharpoonup v^{2^*-1} ~~ \mbox{in} ~~ L^{\frac{2^*}{2^*-1}}(\Omega).
$$
Consequently,
\[
\int_{B_r(0)} {u_n}^{p_n-1}\psi d x \to \int_{B_r(0)} v^{2^* - 1}\psi d x, ~~ \forall \psi \in H^1_0(B_r(0), \re).
\]
Therefore, $v \in H^1_0(B_r(0), \re) \backslash \{0\}$ is a solution of
\[
\begin{cases}
-\lap u + \kappa u = u^{2^* - 1}, \mbox{ in } B_r(0),\\
u = 0, \mbox{ on } \partial B_r(0).
\end{cases}
\]
By Pohozaev's identity, $v \equiv 0$ in $B_r(0)$, and so,
\begin{equation}
\label{phi5}
u_n \tofraco 0 \mbox{ in } H^1_0(B_r(0), \re).
\end{equation}
By definition of $t_n$, we have
\begin{align*}
\int_{B_r(0)} |A(y_n) - A(z + y_n)|^2 |u_n|^{2} & d x + \int_{B_r(0)} |\grad u_n|^2 + \kappa |u_n|^2 d x = \\
                                         = & \int_{\Omega} |A(x) - A(y_n)|^2 |u_n(x - y_n)|^{2} d x  + \\
                                         & + \int_{\Omega} [|\grad u_n(x - y_n)|^2 + \kappa |u_n(x - y_n)|^{2}] d x \\
                                         = & \int_{\Omega} (|\grad_A (e^{i\tau_y}u_n(x - y_n))|^2 + \kappa |e^{i\tau_y}u_n(x - y_n)|^2) d x \\
                                         = & \, {t_n}^{p_n - 2}\int_{\Omega} |e^{i\tau_y}u_n(x - y_n)|^{p_n} d x \\
                                         = & \, {t_n}^{p_n - 2}\int_{B_r(0)} |u_n|^{p_n} d x.
\end{align*}
Since  $u_n \in {\mathcal N}_{\kappa, p_n, r}$,  we get
\begin{equation} \label{Z1}
\int_{B_r(0)} |A(y_n) - A(z + y_n)|^2 |u_n|^{2} d x = ({t_n}^{p_n - 2} - 1)\int_{B_r(0)} (|\grad u_n|^2 + \kappa |u_n|^{2}) d x .
\end{equation}
A direct computation shows that there is $\delta^*>0$ such that
\begin{equation}
\int_{B_r(0)}(|\grad u_n|^2 + \kappa {u_n}^2) dx \geq \delta^* ~~ \forall n \in \mathbb{N}.\label{Z2}
\end{equation}
Combining the boundedness of $(u_n)$ with (\ref{Z1}), (\ref{Z2}), \eqref{cb}, \eqref{phi2} and \eqref{phi5}, we deduce that $t_n \to 1 $.
From \eqref{phi5}, Sobolev embeddings and the boundedness of $(t_n)$,  \eqref{phi1} follows. Since this argument can be applied to any subsequence, the result holds.  \fim


\subsection{Estimates involving the barycenter function}

Consider $\beta: {\mathcal M}_{\kappa, p, \Omega} \to \rn$, the barycenter function, defined as
\[
\beta(u) = \frac{\displaystyle\int_{\Omega} x.|u|^{2^*} d x }{\displaystyle\int_{\Omega} |u|^{2^*} d x }~.
\]

Our first results involving the barycenter function is the following

{\prop \label{bar} For fixed $\kappa \ge 0$, there are $\epsilon = \epsilon(\kappa) > 0$ and $p^* = p^*(\kappa) \in (2, 2^*)$ such that, for $p \in \left.[p^*, 2^*\right.)$, $\beta(u) \in \Omega_r^+$, if $u \in {\mathcal M}_{\kappa, p, \Omega}$ and $I_{\kappa, p, \Omega}(u) \le \frac{1}{N} S^{\frac{N}{2}} + \epsilon$.}

\vspace{0.5 cm}

\noindent {\bf Proof.} ~ Fix $\kappa \ge 0$. By \eqref{cb}, for $p$ close enough to $2^*$, the set
$$
\left\{ u \in {\mathcal M}_{\kappa, p, \Omega} ; I_{\kappa, p, \Omega} (u) \le \frac{1}{N}S^{\frac{N}{2}} + \epsilon \right\}
$$
is non-empty. Suppose, by contradiction, that the result is false. Thus, there are sequences $(p_n)_n, (\epsilon _n)_n$, with $p_n \in (2, 2^*), p_n \to 2^*$ and $\epsilon_n > 0, \epsilon _n \to 0$, and $u_n \in {\mathcal M}_{\kappa, p_n, \Omega}$,  such that
\begin{equation}
\label{bar1}
I_{\kappa, p_n, \Omega} (u_n) \le \frac{1}{N}S^{\frac{N}{2}}+\epsilon_n \mbox{ and } \beta(u_n) \notin \Omega^+_r.
\end{equation}
On the other hand, \eqref{cb} gives
$$
\liminf_{n \to \infty} I_{\kappa, p_n, \Omega} (u_n) \ge \lim_{n \to \infty} b_{\kappa, p_n, \Omega} = \frac{1}{N}S^{\frac{N}{2}}.
$$
Hence, the last two inequalities lead to
\begin{equation}
\label{bar2}
\lim_{n \to \infty} I_{\kappa, p_n, \Omega} (u_n) = \frac{1}{N}S^{\frac{N}{2}}.
\end{equation}
Since $u_n \in {\mathcal M}_{\kappa, p_n, \Omega}$ and $\int_{\Omega} (|\grad_A u_n|^2 + \kappa |u_n|^2) d x = \int_{\Omega} |u_n|^{p_n} d x $, we know that
\[
I_{\kappa, p_n, \Omega} (u_n) = \left(\frac{1}{2} - \frac{1}{p_n}\right) \int_{\Omega} (|\grad_A u_n|^2 + \kappa |u_n|^2) d x
\]
and  by \eqref{bar2},
$$
\lim_{n \to \infty} \int_{\Omega} (|\grad_A u_n|^2 + \kappa |u_n|^2) d x = S^{\frac{N}{2}}.
$$
The above limit yields
\[
\lim_{n \to \infty}\frac{\displaystyle\int_{\Omega} (|\grad_A u_n|^2 + \kappa |u_n|^2) d x }{\left(\displaystyle\int_{\Omega} |u_n|^{p_n} d x \right)^{\frac{2}{p_n}}} = \lim_{n \to \infty}\left( \int_{\Omega} (|\grad_A u_n|^2 + \kappa |u_n|^2) d x \right)^{1 - \frac{2}{p_n}} = S.
\]
Using the diamagnetic inequality and the last limit, we get
\begin{equation}
\label{bar3}
\limsup_{n \to \infty} \frac{\displaystyle\int_{\Omega} (|\grad |u_n||^2 + \kappa |u_n|^2) d x }{\left(\displaystyle\int_{\Omega} |u_n|^{p_n} d x \right)^{\frac{2}{p_n}}} \le \lim_{n \to \infty} \frac{\displaystyle\int_{\Omega} (|\grad_A u_n|^2 + \kappa |u_n|^2) d x }{\left(\displaystyle\int_{\Omega} |u_n|^{p_n} d x \right)^{\frac{2}{p_n}}} = S
\end{equation}
The limit \eqref{bar3} implies that, for $\delta_1 > 0$ to be chosen later, there is $n_{1} \in \mathbb N$ such that for $n \ge n_{1}$,
\begin{equation}
\label{bar4}
\frac{\displaystyle\int_{\Omega} (|\grad |u_n||^2 + \kappa |u_n|^2) d x }{\left(\displaystyle\int_{\Omega} |u_n|^{p_n} d x \right)^{\frac{2}{p_n}}} \le S + \delta_1.
\end{equation}
Arguing as in \eqref{bkpd3}, for $\delta_2 > 0$  to be also chosen later, there is $n_{2} \in \mathbb N$ such that for $n \ge n_{2}$,
\begin{equation}
\label{bar5}
\frac{\displaystyle\int_{\Omega} (|\grad |u_n||^2 + \kappa |u_n|^2) d x }{\left(\displaystyle\int_{\Omega} |u_n|^{2^*} d x \right)^{\frac{2}{2^*}}} \le \frac{\displaystyle\int_{\Omega} (|\grad |u_n||^2 + \kappa |u_n|^2) d x }{\left(\displaystyle\int_{\Omega} |u_n|^{p_n} d x \right)^{\frac{2}{p_n}}} + \delta_2.
\end{equation}
From \eqref{bar4} and \eqref{bar5}, for $n \ge \max\limits_{j = 1, 2} n_{j}$, we have
\begin{equation}
\label{bar6}
S \le \frac{\displaystyle\int_{\Omega} (|\grad |u_n||^2 + \kappa |u_n|^2) \dd x }{\left(\displaystyle\int_{\Omega} |u_n|^{2^*} \dd x \right)^{\frac{2}{2^*}}} \le S + \delta_1 + \delta_2.
\end{equation}
We claim that  there is $\eta > 0$ such that if $v \in H^1_0(\Omega)$ satisfies
\begin{equation}\label{claim}
\frac{\displaystyle\int_{\Omega} (|\grad v|^2 + \kappa v^2) \dd x }{\left(\displaystyle\int_{\Omega} |v|^{2^*} \dd x \right)^{\frac{2}{2^*}}} \le S + \eta,
\end{equation}
then  $\beta(v) \in \Omega^+_r.$ Effectively,  suppose by contradiction that (\ref{claim}) does not hold. Thus, there are $(v_n) \subset H^1_0(\Omega, \re)$ and $\eta_n \to 0$ such that
\[
\frac{\displaystyle\int_{\Omega} (|\grad v_n|^2 + \kappa |v_n|^2) d x }{\left(\displaystyle\int_{\Omega} |v_n|^{2^*} d x \right)^{\frac{2}{2^*}}} \le S + \eta_n, \mbox{ with } \beta(v_n) \notin \Omega^+_r.
\]
Let $w_n := {v_n}/{|v_n|_{2^*, \Omega}}$. Thus, $(w_n) \subset H^1_0(\Omega)$ is a bounded sequence. Hence,  there are $u \in H^1_0(\Omega)$ and finite positive measures $\mu, \nu \in \mathcal M (\rn)$ verifying, for some subsequence,
\[
\begin{cases}
|w_n| \tofraco u \mbox{ in } D^{1,2}(\rn), \mbox{ as } n \to \infty,\\
|\grad w_n - \grad u|^2 \tofraco \mu \mbox{ in } \mathcal M (\rn), \mbox{ as } n \to \infty,\\
|w_n - u|^{2^*} \tofraco \nu \mbox{ in } \mathcal M (\rn), \mbox{ as } n \to \infty,\\
w_n(x) \to u(x) \mbox{ almost everywhere } \Omega, \mbox{ as } n \to \infty,
\end{cases}
\]
where we made the extension by zero outside of $\Omega$. By Concentration-compactness lemma,
\[
S = |\grad u|_{2}^2 + \|\mu\|_{\mathcal M (\rn)}, \,\, 1 = |u|_{2^*}^{2^*} + \|\nu\|_{\mathcal M (\rn)}, \,\, \|\nu\|_{\mathcal M (\rn)}^{\frac{2}{2^*}} \le S^{-1}\|\mu\|_{\mathcal M (\rn)}.
\]
Employing the arguments in \cite{Willem}, $\nu$ and $\mu$ are concentrated at $y \in \overline{\Omega}$ and satisfy $\|\nu\|_{\mathcal M (\rn)}^{\frac{2}{2^*}} = S^{-1}\|\mu\|_{\mathcal M (\rn)}$. Let $\Gamma:\rn \to \rn$ and $\Upsilon:\rn \to \re$ be continuous functions with compact support such that in a neighborhood of $\overline{\Omega}$, $\Gamma = Id_{\rn}$ and $\Upsilon = 1$. Using these functions, we derive
$$
\beta(v_n) = \beta(w_n) = \frac{\displaystyle\int_{\Omega} x.|w_n|^{2^*} d x }{\displaystyle\int_{\Omega} |w_n|^{2^*} d x } = \frac{\displaystyle\int_{\rn} \Gamma(x) |w_n|^{2^*} d x }{\displaystyle\int_{\rn} \Upsilon(x) |w_n|^{2^*} d x }
$$
Hence,
$$
\beta(v_n) \to \frac{\displaystyle\int_{\{y\}} \Gamma(x) d \nu}{\displaystyle\int_{\{y\}} \Upsilon(x) d \nu} = \frac{\nu (y) \Gamma(y)}{\nu(y) \Upsilon(y)} = y \in \overline\Omega,
$$
contradicting the fact that $\beta(v_n) \notin \Omega$. Hence, the (\ref{claim}) holds.  For $\eta$ given by (\ref{claim}), take in \eqref{bar6}, $\delta_1 = \delta_2 = \displaystyle\frac{\eta}{2}$. Observing that  $\beta(|u_n|) = \beta(u_n)$, we have,
\[
\beta(u_n) \in \Omega^+_r,
\]
which contradicts (\ref{bar1}) and the proof is complete. \fim

\vspace{0.5 cm}

For any $\kappa \ge 0$ fixed, consider  $\epsilon = \epsilon(\kappa) > 0$ given by Proposition \ref{bar}. Define
\begin{equation}
\label{e*}
\epsilon ^* = \epsilon^* (\kappa) = \frac{1}{N}S^{\frac{N}{2}} + \epsilon
\end{equation}
and the set
\[
{\mathcal M}_{\kappa, p, \Omega}^{\epsilon^*} := \{u \in  {\mathcal M}_{\kappa, p, \Omega} ; I_{\kappa, p, \Omega}(u) \le \epsilon^*\}.
\]

\begin{corollary} \label{pb} For fixed $\kappa \ge 0$, there is $\overline p(\kappa) \in (2, 2^*)$ such that, for each $p \in \left.[ \overline p(\kappa), 2^*\right.)$,
\[
\Phi_{\kappa, p}(\Omega^-_r) \subset {\mathcal M}^{\epsilon^*}_{\kappa, p, \Omega}, \quad \beta({\mathcal M}_{\kappa, p, \Omega} ^{\epsilon^*}) \subset \Omega^+_r.
\]
\end{corollary}

\noindent {\bf Proof.} The proof follows immediately from Lemma \ref{phi} and Proposition \ref{bar}. \fim

\subsection{Proofs of Theorem \ref{t1} and Corollary \ref{c1}}

We are now ready to conclude the proof of Theorem \ref{t1}.  The key ingredient is the verification of Theorem \ref{resabs}.  To this end,  fix  $\kappa \ge 0$,  and  take  $p \in [\left.\overline p, 2^*)\right.$,  for $\overline p = \overline p(\kappa)$ given by Lemma \ref{pb}. Let $\mathcal K$ be the set of critical points of $I_{\kappa, p, \Omega}$. Suppose that $\mathcal K$ is discrete.
We begin observing that condition $(i)$  is a consequence of the definition of $I_{\kappa, p, \Omega}$, for $\Psi$ given by
$\Psi(u) = \frac{1}{p}\int_\Omega|u|^p\,dx$.
Using that the Hessian form of $I_{\kappa, p, \Omega}$ at $u$ is given by
 \[
H_{I_{\kappa, p, \Omega}}(u)(v, w) = \langle v, w \rangle_{E} - (p - 1) \int_{\Omega} |u|^{p-2} \operatorname{Re} (w\overline v) \dd x, \quad \forall v, w \in E,
\]
we have that $H_{I_{\kappa, p, \Omega}}(u)$ is a bounded symmetric bilinear form, for every $u \in E$. The Riesz representation produces a self-adjoint operator $L(u): E \to E$ such that $H_{I_{\kappa, p, \Omega}}(u)(v, v) = \langle L(u)v, v\rangle_E$.  This and Proposition \ref{ps} imply that condition $(ii)$ holds.
By Proposition \ref{nehari}, the Nehari manifold $\mathcal{M}_{\kappa,p,\Omega}$ is homeomorphic to the unit sphere of $E$, which implies $(iii)$. Consider $\epsilon^*$ given by \eqref{e*}. We can clearly assume that  $\epsilon^*$ is a regular level of $I_{\kappa, p, \Omega}$. By Corolary \ref{pb}, for $p \in [\left.\overline p, 2^*)\right.$ the maps  $\Phi_{\kappa, p} : \Omega^-_r \to {\mathcal M}^{\epsilon^*}_{\kappa, p, \Omega}$ and $\beta : {\mathcal M}^{\epsilon^*}_{\kappa, p, \Omega} \to \Omega_r^+$ are continuous and satisfy $\beta \circ \Phi_{\kappa, p} = Id_{\Omega_r^-}$, where, by construction,  $\Omega^+_r, \Omega^+_r$ are homotopically equivalent to $\Omega$. We conclude that $(iv)$ holds. Consequently, by  Theorem \ref{resabs}, we have
\begin{equation*}
\sum_{u \in \, \mathcal C_1} i_t (u) = t \mathcal P_t(\Omega) + t \mathcal Q(t) + (1+t)\mathcal Q_1(t)
\end{equation*}
and
\begin{equation*}
\sum_{u \in \, \mathcal C_2} i_t (u) = t^2 [\mathcal P_t(\Omega) + \mathcal Q(t) -1 ] + (1+t)\mathcal Q_2(t),
\end{equation*}
where, for $\delta \in (0, \delta)$, $\delta > 0$ given by Proposition \ref{nehari},
\[
\mathcal C_1 := \{ u \in \mathcal K; \delta < I_{\kappa, p, \Omega}(u) \le \epsilon^* \}, \quad \mathcal C_2 := \{ u \in \mathcal K; \epsilon^* < I_{\kappa, p, \Omega}(u)\}.
\]
Thus
\[
\sum_{u \in \mathcal K} i_t (u) = t \mathcal P_t(\Omega) + t^2 [\mathcal P_t(\Omega) - 1 ] + \mathcal Q_3(t),
\]
where $\mathcal Q _3$ is a polynomial with non-negative coefficients. The proof of Theorem \ref{t1} is complete. In order to prove Corollary \ref{c1}, suppose that every critical point of $I_{\kappa, p, \Omega}$ is non-degenerate. By general Morse theory,
\[
i(u) = t^{m(u)}, \mbox{ for all } u \in \mathcal K,
\]
and the result follows from Theorem \ref{t1}.
\cqd


\end{document}